American University of Armenia

College of Science and Engineering

Numerical Analysis Project

# The numerical solution of Navier hyperbolic equations: Shock wave propagation

Arakelyan Erik, Jilavyan Narek, Serobyan Aram

Fall, 2015

# Abstract

In the following paper we will consider Navier-Stokes problem and it's interpretation by hyperbolic waves, focusing on wave propagation. We will begin with solution for linear waves, then present problem for non-linear waves. Later we will derive for numerical solution using PDE's. Also we will design a Matlab program to solve and simulate wave propagation.

**Keywords[1]:** Navier-Stokes equations, Fluid dynamics, non-linear waves, linear waves, hydrodynamics, hyperbolic waves, wave simulation, Ordinary differential equations, Partial differential equations, shock wave propagation, wave equation, string vibration, signal propagation, Explicit Finite difference approximation.

---

[1] And fancy stuff :)

# Contents

# List of Figures

# Chapter 1

# Introduction

## 1.1 Problem Setting and Description

Navier-Stokes equations are designed in order to simulate the flow of various fluids in numerous disciplines across physics (electromagnetism, acoustics, elasticity). Navier- Stokes equations were formed by combining Newton's second law and fluid motion. Also assumption of stress being the sum of diffusing viscous term and pressure term was taken to describe viscous flow. Without Navier-Stokes equations working with weather model, ocean currents, water flow in a pipe, air flow around a wing would be extremely hard. Also these equations are widely used in designing airplanes and cars, studying blood flow, designing power stations, analysing pollution and so on.

History of the Navier-Stokes equations dates back to 1822. Claude-Louis Navier mentioned in his memoir the equations for homogeneous incompressible fluids with regards to molecular arguments. Later, in 1829 Poisson derived equations for compressible fluids. Consequentially, George Gabriel Stokes formulated the Navier-Stokes equations in his paper (1847) as we know it today.

Up to this day there is no proof for existence and uniqueness of Navier-Stokes equations' solution in three dimensional space. Moreover, it is considered one of the seven most important unsolved problems currently. Clay Mathematical Institute offers a significant prize of one million dollars to the person who will prove this or bring a counter-example. Although a number of advancements have been made in this sphere, yet it remains mainly a "dark forest" for the researches. As you will see in this paper, there is a perfectly working solution for both one dimensional case and two dimensional case.

This paper illustrates Navier-Stokes hyperbolic wave equations, its solutions in one dimensional case, and extending that solution for two dimensional cases. The solution for one dimensional case is rather simple, in contrast the two-dimensional case requires more work. In following paper we focused on solution by using partial differential equations. Explicit finite difference approximation was used for that purpose. In the end we bring wave propagation examples, with both two dimensional and three dimensional figures.

## 1.1.1 The structure of the project paper

This paper consists of six chapters. It begins with introductory chapter, which you are reading right now. Then it proceeds to chapter two, where slight theoretical background is given about hyperbolic wave equations. There is also an analytical solution for one dimensional case of the problem. Next, in chapter three we solve numerically two dimensional case. In chapter four, you can find the solution for some specific problems, and learn about the steps to be taken in the future. The following chapter contains two appendices. First is the glossary of the unknown terms, and in second appendix you can find our Matlab codes. Finally, last, the sixth chapter contains reference to all materials we have used.

# Chapter 2

# Theoretical Background

## 2.1 Hyperbolic Waves

In our paper we are interested in hydrodynamic equations which are similar with signal-propagation equations. That specific kind of equations are called hyperbolic equations. Also hyperbolic equations have much broader use, and can be met elsewhere, we will concentrate on them used in hydrodynamics only.
The simplest hyperbolic equations are given by the following equation

$$\varphi_{tt} = c_0^2 \nabla^2 \varphi \quad (0)$$

where φ is a function φ(x,t), $c_0$ is constant velocity, and $\nabla^2 = \Delta f = \sum_{i=0}^{n} \frac{\partial^2 f}{\partial x_i^2}$. In (2.2) we will derive the following formula, representing hyperbolic waves in one dimensional case.

$$\varphi_t + c_0 \cdot \varphi_x = 0 \quad (1)$$

**(1)** is the common equation for linear waves. Here φ(x, t) is a function which describes time depending shifting of the function φ along x which has a velocity $c_0$. That is it describes the propagation of the wave.

Let's understand difference between linear waves and non-linear. In linear waves the velocity at any point is the same. That is, the total velocity of the flow is not effected from internal particle velocities (which we will see are the case in non-linear). So if we look at *Figure 2.1,* we'll see that the wave repeats itself and it propagates without any change at any moment in time. If we would pick any point and examine it during time, the overall velocity in that point will always be the same. For example, imagine that you are in the middle of the ocean. Waves are born and start their propagation. From that point waves are linear as they acquire some height and just propagate without changing their shape.

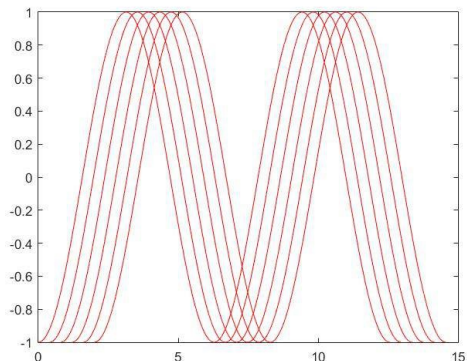

*Figure 2.1*

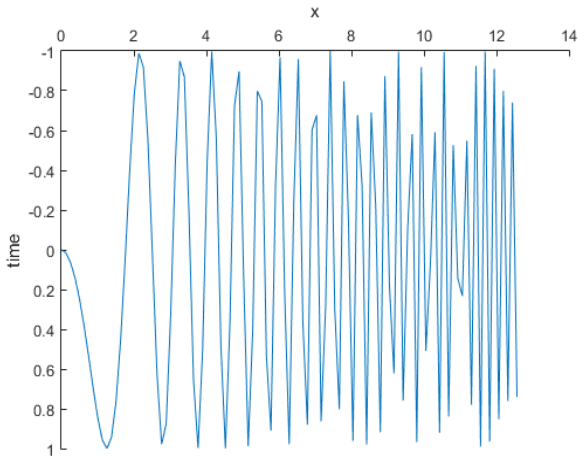

*Figure 2.2*

What about nonlinear waves, we will rarely observe too similar curvatures at two different moments. Particles in one position have higher velocity than in others, thus during the time we reach a point when some particle x, reaches particle y, which was in front before, then passes it. We will have then a picture of crushing wave.
For instance, imagine you are on the beach. You see the waves coming, then crushing on the sand. Such waves are not linear, as they change their shape over time.
In *Figure 2.2* you can see nonlinear wave. As in computers for crushing wave one x would require three y values, it is not display, instead such image is resulted.

## 2.2 Linear Waves

In 2.1 we presented the general equation for hyperbolic waves. Now let's consider it for one dimensional case. First let's rewrite **(0)** for one dimension.

$$\varphi_{tt} = c_0^2 \varphi_{xx}$$

Then rewrite it in the following form.

$$\varphi_{tt} - c_0^2 \varphi_{xx} = 0$$

Now, let's perform variable change with variables α and β such that x = $c_0$t + α, x = β - $c_0$t. Let's take φ(x) to be equal the sum of two functions f(α) and g(β). Here f(α) denotes a function of a wave propagating to the right, and g(β) to the left. As can be seen, β = $c_0$t + x, α = $c_0$t - x.

$$\varphi = f(\alpha) + g(\beta) = f(x - c_0 t) + g(x + c_0 t) \quad \textbf{(2)}$$

From **(2)** we can see that φ(x, t) is represented as a summation of 2 independent functions, thus $\varphi_{\alpha\beta}$ =0. (Witham, 1974)

The **(2)** is general solution for one dimensional case however we face a problem: we have two waves. Nevertheless our problem has a simple solution as there is a way to get rid of one of the waves.

$$\varphi_{\alpha\beta} = \frac{\partial^2 \varphi}{\partial \alpha \partial \beta} = \frac{\partial}{\partial \alpha}\frac{\partial}{\partial \beta}\varphi = \left(\frac{\partial}{\partial t} - c_0 \frac{\partial}{\partial x}\right)\left(\frac{\partial}{\partial t} + c_0 \frac{\partial}{\partial x}\right)\varphi = 0$$

As we get that the derivative of the function with respect α and then β equals to 0, it is enough keeping only one of the equations in the parenthesis. Thus we get alternative form of **(1)** i.e.

$$\frac{\partial \varphi}{\partial t} + c_0 \frac{\partial \varphi}{\partial x} = 0 \quad \textbf{(3)}$$

This is the first order wave equation. In that equation φ = φ(x, t), t, x are our lovely time-axis and x-axis, and $c_0$ is the constant velocity. For the sake of making it more durable for ordinary person let's examine the flow of the river. When discussing fluid systems one can describe the whole system as a sum of individual particles. An interesting phenomenon occurs in this kind of systems. Although the velocity of individual particles might change drasticly in the river over time, the overall velocity of the river remains the same. If the speed of the river would not remain the same, then the wide and narrow parts would not remain in same place.

Now, let's transform **(3)** to the system of ordinary differential equations. For that we substitute $|c_0|$ (taking the absolute value because the velocity cannot be negative) by dx/dt. Obviously it is constant as it is the absolute value of the derivative of x - $c_0$t with respect to t.
The following system is derived.

$$\begin{cases} \frac{\partial \varphi}{\partial t} + \frac{dx}{dt}\frac{\partial \varphi}{\partial x} = 0 \\ \frac{dx}{dt} = |c_0| \end{cases} \Rightarrow \begin{cases} \frac{d\varphi}{dt} = 0 \\ \frac{dx}{dt} = |c_0| \end{cases} \quad \textbf{(4)}$$

From **(4)** we conclude that velocity of φ is changing with constant velocity $c_0$. Meaning that wave propagates from left to right with constant velocity $c_0$.

- 7 -

## 2.3 Non-linear Waves

In order to make things more interesting, let's make them a little more complicated. For that very purpose we substitute the constant velocity change $c_0$ by non-constant function $\varphi$. Thus we get the equation:

$$\frac{\partial \varphi}{\partial t} + \zeta \frac{\partial \varphi}{\partial x} = 0 \quad \textbf{(5)}$$

Here $\zeta$ can be a function depending on time, x , or $\varphi$ itself. Now, not only velocity of individual particles changes, but also the change of the velocity is described by non-constant function. To this date no one is known to be close to solving this analytically, which is not surprising as it is incredibly hard, if not impossible to solve. Nevertheless it is possible to come up with numerical solution, which we did in next chapter. Also there are two possible approaches we can choose for solving **(5)**. First approach involves describing **(5)** using Ordinary differential equations and giving numerical solution for them. The other method is solving PDE's numerically using substitutions.

Solution using ODE's:
Deriving ODE's from **(5)** doesn't raise any difficulty, we just plug $c_0 = \zeta(x,t)$ (or other function satisfying **(4)**) and get the following result.

$$\begin{cases} \dfrac{d\varphi}{dt} = 0 \\ \dfrac{dx}{dt} = \zeta \end{cases} \quad \textbf{(6)}$$

Equation **(6)** correctly represent wave prorogation for waves with amplitudes smaller than one, however it's rather ambitious representation for waves with higher amplitudes. In this paper we didn't concentrated on solving **(6)**.

Solution using approximation for PDE's:

The other way for solving **(5)** is choosing small enough $\Delta t$ in order to be able to neglect the limit for approximation of the PDE. We know that

$$\frac{\partial \varphi}{\partial t} = \lim_{\Delta t \to 0} \frac{\varphi(t+\Delta t, x) - \varphi(t,x)}{\Delta t} \approx \frac{\varphi(t+\Delta t, x) - \varphi(t,x)}{\Delta t}$$

and

$$\frac{\partial \varphi}{\partial x} = \lim_{\Delta x \to 0} \frac{\varphi(t, x+\Delta x) - \varphi(t, x)}{\Delta x} \approx \frac{\varphi(t, x+\Delta x) - \varphi(t, x)}{\Delta x}$$

the last one also can be represented as

$$\frac{\partial \varphi}{\partial x} = \lim_{\Delta x \to 0} \frac{\varphi(t, x+\Delta x) - \varphi(t, x-\Delta x)}{2\Delta x} \approx \frac{\varphi(t, x+\Delta x) - \varphi(t, x-\Delta x)}{2\Delta x}$$

Thus, to not be biased to the right or left, we write it in following way

$$\lim_{\Delta x \to 0} \frac{\varphi(t, x+\Delta x) - \varphi(t, x-\Delta x)}{\Delta 2x} \approx \frac{\varphi(t, x+\Delta x) - \varphi(t, x-\Delta x)}{\Delta 2x}$$

plugging this into **(5)** we get

$$\frac{\varphi(t+\Delta t, x) - \varphi(t, x)}{\Delta t} - \zeta \frac{\varphi(t, x+\Delta x) - \varphi(t, x-\Delta x)}{2\Delta x} = 0 \quad \textbf{(7)}$$

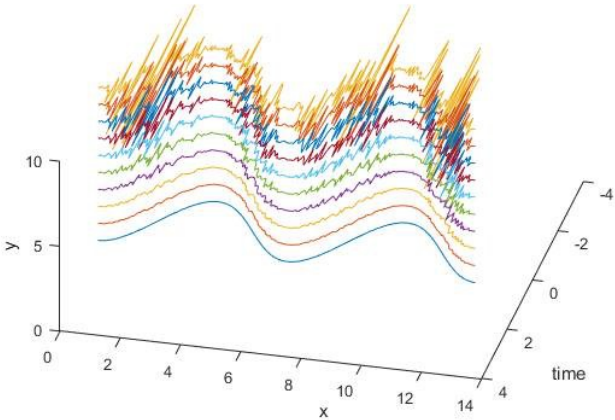

physical meaning of this equation is the propagation of non-linear wave, however in some parts the amount of error can cause grainy (non-smooth) surface in the graph. As can be seen from *Figure 2.3* the wave propagating over time starts to destabilize due to computational-error accumulation.

*Figure 2.3*

# Chapter 3

# Numerical Solution

## 3.1 Numerical Background and step by step algorithms

For 1D case no special calculating method must be introduced as it is solved analytically by using, D'Alembert's solution that is presented in **2.2**.

$$\varphi(x) = f(x - c_0 \cdot t)$$

However, when we get to solve the PDE for the non-linear wave propagation, there is a methodology to be presented. As we already presented in paragraph **2.3** when having a hyperbolic PDE of type

$$\lim_{\Delta t \to 0} \frac{\varphi(t+\Delta t, x) - \varphi(t, x)}{\Delta t} - \zeta \lim_{\Delta x \to 0} \frac{\varphi(t, x+\Delta x) - \varphi(t, x-\Delta x)}{2\Delta x} = 0$$

we can approximate it by neglecting the limit and taking a small Δt and Δx. As we do we get the approximation

$$\frac{\varphi(t+\Delta t, x) - \varphi(t, x)}{\Delta t} - \zeta \frac{\varphi(t, x+\Delta x) - \varphi(t, x-\Delta x)}{2\Delta x} = 0$$

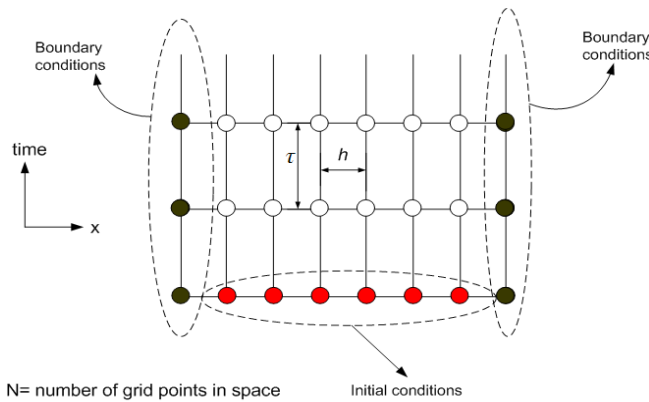

After getting this equation lets divide our coordinate system into an n by m grid where where the difference between each to $x_i$ neighbors is Δx and $t_i$ neighbors Δt.

Now having this grid we can rewrite the equation we got in the following form.

$$\frac{\varphi(t_{i+1}, x) - \varphi(t, x)}{\Delta t} - \zeta \frac{\varphi(t, x_{i+1}) - \varphi(t, x_{i-1})}{2\Delta x} = 0$$

From here, given the initial condition (row φ(0,x)), we want to find the values of each row to reconstruct the wave at each point of time.

That is we need to find φ(t$_{i+1}$,x). Lets get it from the equation above.

$$\varphi(t_{i+1}, x) = \Delta t \cdot \zeta \frac{\varphi(t, x_{i+1}) - \varphi(t, x_{i-1})}{2 \Delta x} + \varphi(t, x)$$

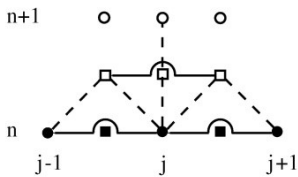

Here is the numerical solution for the The PDE that we intended to solve. As one can notice this scheme can have a problem at the boundaries. Numerous methods have been developed to approximate and solve that problems. The famous schemes are Lax-Wanderof, leap-frog etc etc.

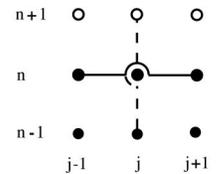

## 3.2 Numerical Solution of Navier-Stokes for linear waves

The numerical solution in the this case involves just using the analytical solution and simulating time in order to get the wave propagation. One just needs to assign $\varphi(x) = f(x - c_0 \cdot t)$ and correctly simulate the propagation of wave with the constant speed $c_0$. For this simulation one can use the horizontal translation of the function by $c_0 t$.

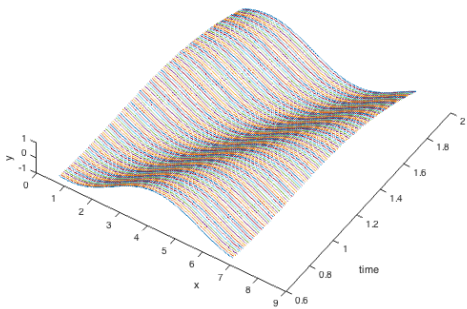

As one can see from *Figure 3.2.1* the wave propgates with constant speed. The periodic for of the wave is the obvious assertion for the abovementioned words. Also no unsmoothness or error cumulation is theoretically possible in this case, due to the fact that the solution for the equation is derived analytically.

*Figure 3.2.1*

## 3.3 Numerical Solution of Navier-Stokes for non-linear waves

In this section we will speak about main aim of this paper, the solution of two dimensional that is non-linear waves. As was mentioned in the previous chapter we will derive solution by PDEs. First let's understand what we have and what we want to get. We have the function at the initial moment (t=0), and we want to be able to trace how it would behave during some given period. Also we have function ζ which describes the velocity change for the individual particle at the given time moment. **(7)** shows relation we had derived.

For our purpose we will have to approximate the initial function at each next stage by result of the previous step. To make it more clear let's examine *figure 3.3.1*. In the given example we have propagation of function $\sin(x^2)$ at two time moments. Also we are looking at the range [0, 4π]. As we see the waves are similar however they slightly differ at the endpoints. The reason for alternative propagation is mainly because of ζ. In the figure it is defined by x+t, which forces the wave to change its shape.

The change in the wave is seen much better in *figure 3.3.2* which demonstrates the wave in 16 different time moments. Two different view are presented to show the difference in the shapes.

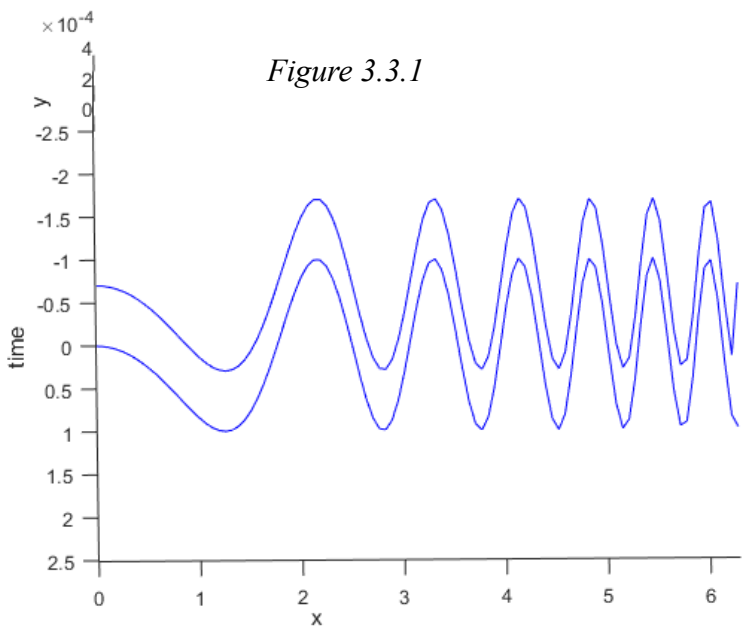

*Figure 3.3.1*

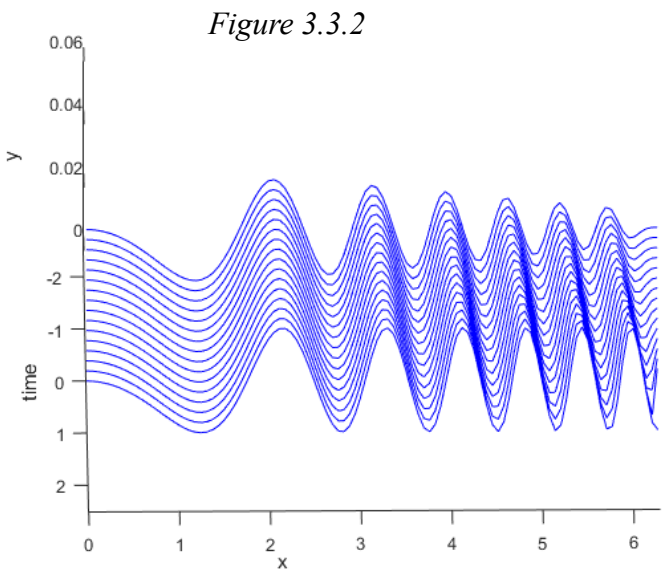

*Figure 3.3.2*

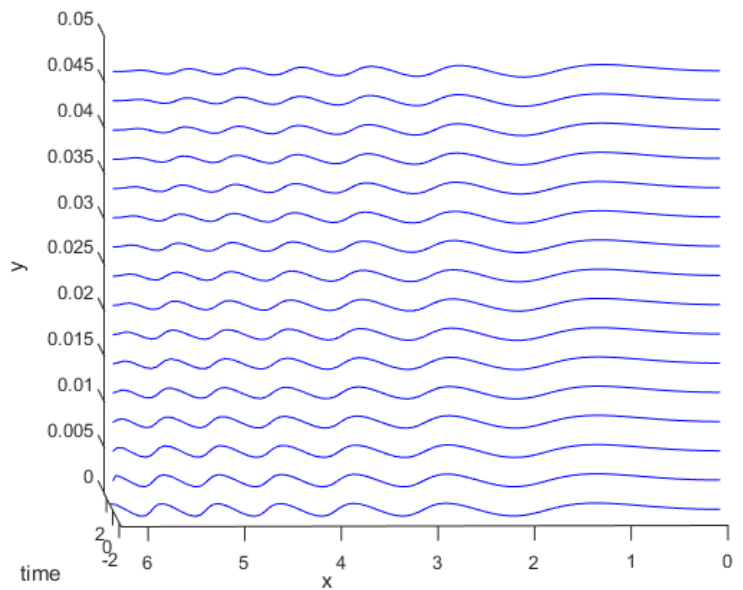

There are several methods for such approximation. In this paper we chose two methods for doing it. There is the same idea behind both methods however the construction slightly differs for each other step. Beginning of steps is similar. We take two dimensional array. There is t rows in the array, where t stands for time steps ( Δt). Also width is the length of the desired interval divided into some n. n is predefined usually large number. The amount of points taken for constructing the next function is equal to the width of the array. This much about basic idea of the algorithm.

Question arises, how we construct function for next step. We had **(7)** which we can rewrite in the following form.

$$\frac{\varphi(t+\Delta t, x)-\varphi(t,x)}{\Delta t} = \zeta \frac{\varphi(t,x+\Delta x)-\varphi(t,x-\Delta x)}{2\Delta x}$$

Our aim is to approximate the function for next time moment. For that purpose we rewrite it again in the following way.

$$\varphi(t+\Delta t,x) = \Delta t \cdot \zeta \frac{\varphi(t,x+\Delta x)-\varphi(t,x-\Delta x)}{2\Delta x}+\varphi(t,x) \quad \textbf{(8)}$$

We acquired a formula for approximating a point of the wave in the next time moment using three points from the wave at that moment. However, as we aim to make our calculations using computer program there is an issue that must be solved. As we are using three points from the previous time moment we either need to declare array twice larger than amount of desired points or find some other way for finding the boundary points. We decided to pick the second way, thus we approximate such points either by the closest point in that moment, or using **(8)** and plugging instead of missing construction point the closest point to it at that time moment.

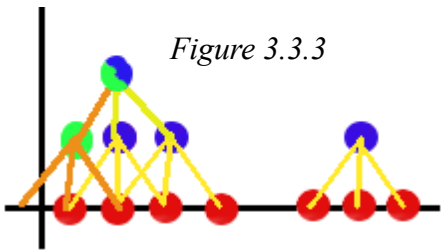

*Figure 3.3.3*

Look at *Figure 3.3.3*. The red points are the initial points, that is points at t = 0. The blue points, are points of our function at t = 0+Δt. The yellow lines connect points which were used for construction of the blue point.

As you see there is a green point, which has an orange line going nowhere. In our code we are solving that issue either by taking the blue point next to it, or by **(8)** with double usage of the point under it. Though there is a larger error at each time moment.

### Our functions:

oneD is a function for 3d plotting of linear waves.

tp1 plots in 3D wave at different moment. It uses the upper mention method, the boundary points it replaces by closes point at that time moment.

tp2 plots in 3D again however it uses twice one of the points ( in *figure 3.3.3* it is the point in the bottom)

tp3 plots 2D, showing wave propagation. It uses same structure as tp2.
tp4 is the 2D version of tp1 which illustrates its propagation.

Doing multiple experiments we found out that structure of tp2 gives a better result than tp1. You can find our results in Chapter 4.

# Chapter 4

# The Results and Conclusion

## 4.1 Results for 1D.

For one dimensional case our results are very stable. First let's examine function sinx with constant speed 1 on [0,4π], given t=10.

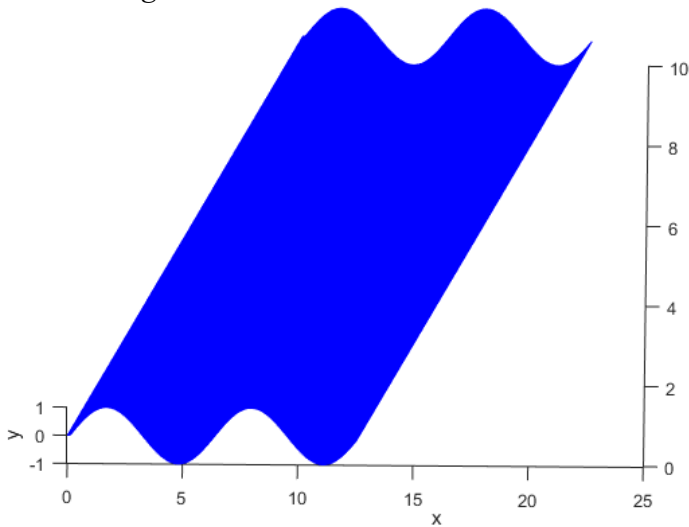

*Figure 4.1.1*

*Figure 4.1.1* is the plot by using OneD function. As we see our wave propagates without any error.

*Figure4.1.2 Figure* is the plot using function tp1, as we see even for linear waves there exist an error caused by approximation of border points.

*Figure4.1.3* Again there is no error in propagation. tp2 is used.

*Figure 4.1.4a,b* this figures illustrates the propagation of wave constructed by tp3 in different time moments in 2D and *Figure4.1.5a,b* illustrates the propagation of wave using tp4. Error accumulation is obvious.

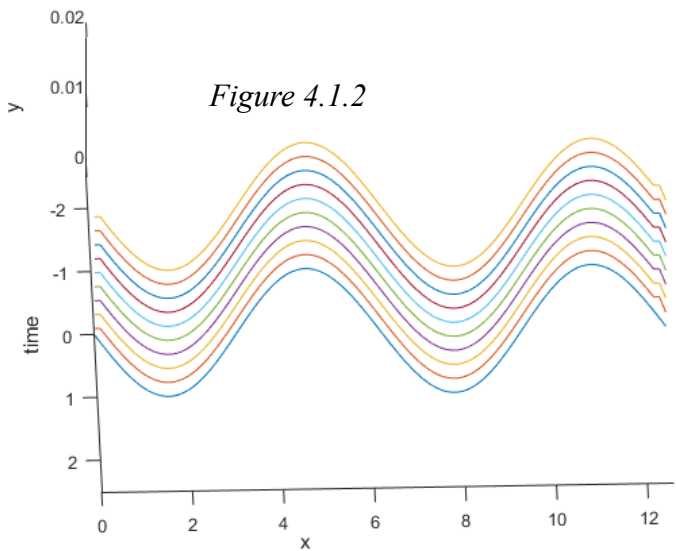

*Figure 4.1.2*

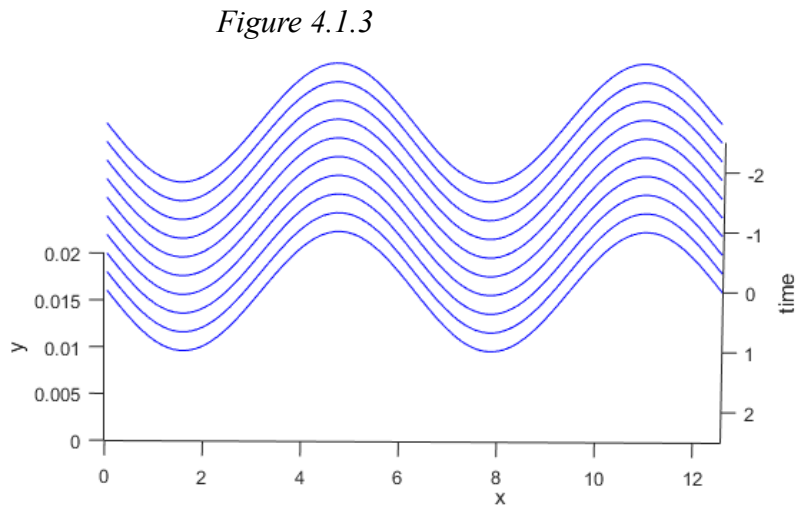

*Figure 4.1.3*

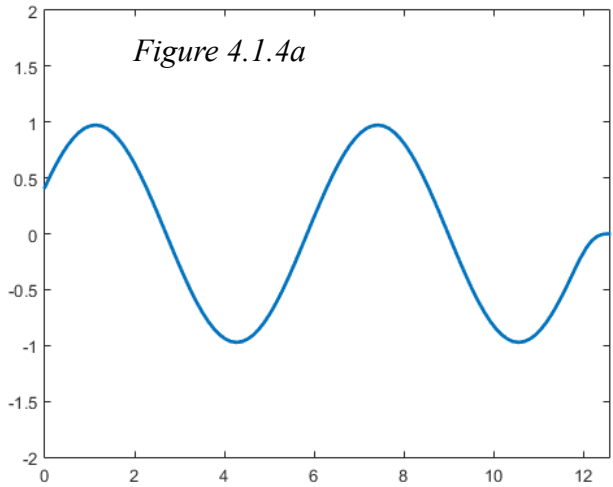
*Figure 4.1.4a*

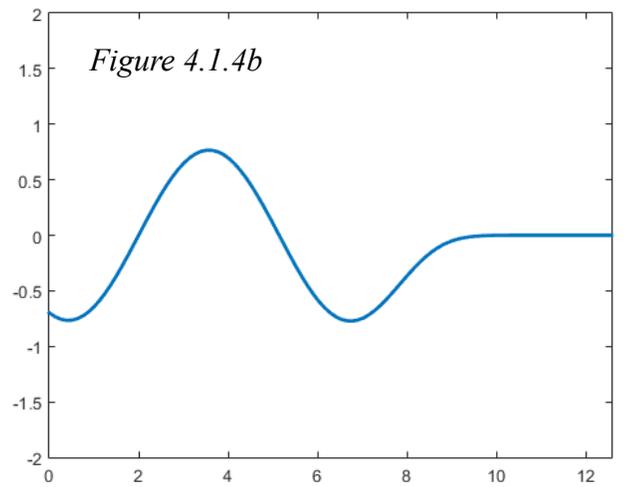
*Figure 4.1.4b*

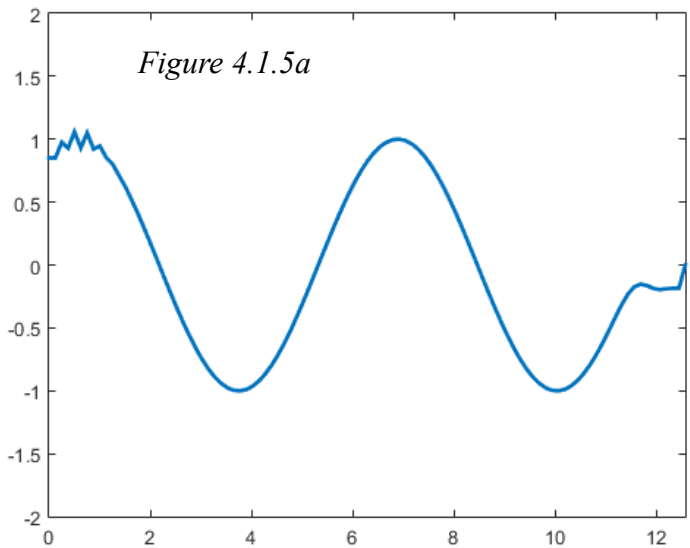
*Figure 4.1.5a*

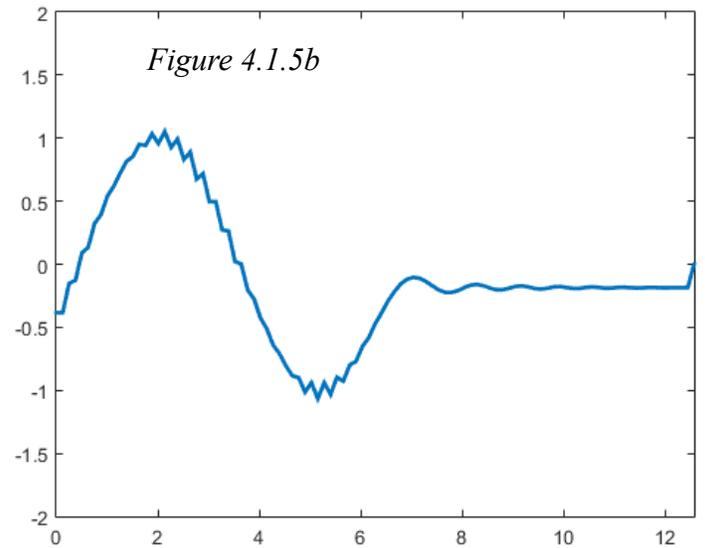
*Figure 4.1.5b*

There is no need for further exploration of linear wave case, thus we will move on to two dimensional case.

## 4.2 Results for 2D.

In Section 4.1 we witnessed that structure of tp1gives not very good results, but let's start this Section by couple examples how bad it is.

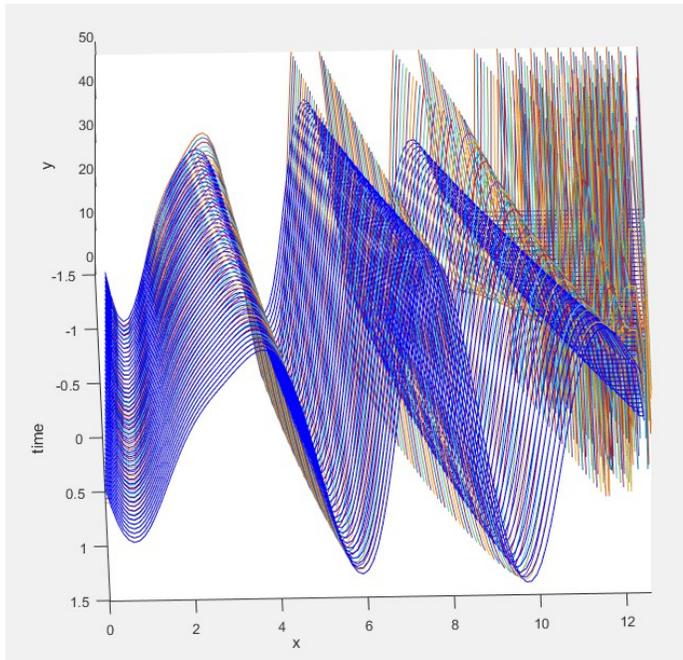

In *figure 4.2.1* blue wave is the wave constructed by tp2. We see that it looks like a normal wave which is narrowing down, however the other wave (with multiple color) behaves very bad and is distorted.

In *figure 4.2.2* and *4.2.3* Illustration of same type of error for tp3 and tp4 is given. I think there is no more point for investigating by both methods and from now on examples will be constructed exclusively by tp3 and tp2.

For all three figures $sin(x^2)$ function is taking with speed depending from both t and x given by x+t. Range is $[0,4\pi]$.

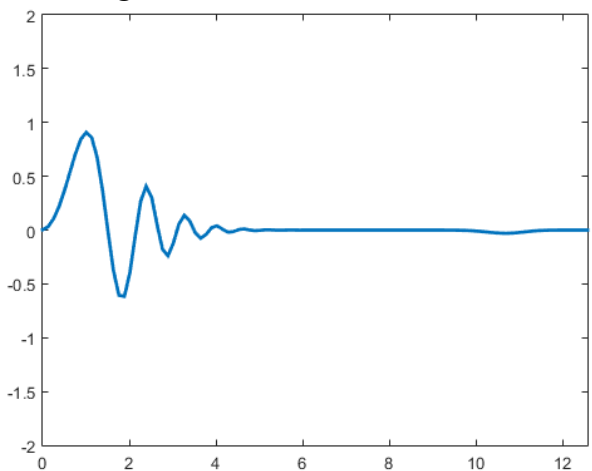

*Figure 4.2.2*

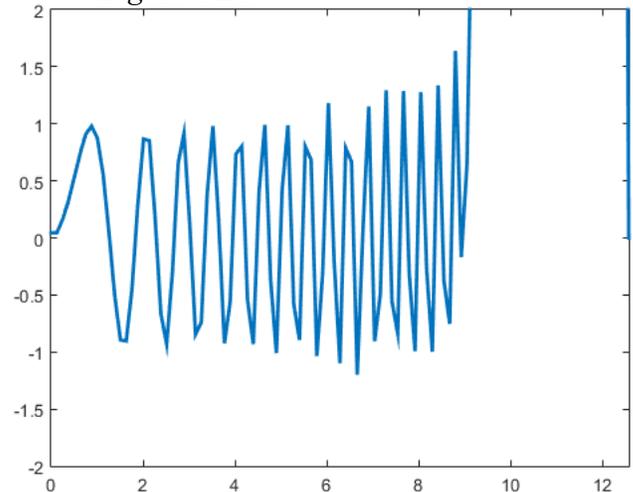

*Figure 4.2.3*

There are several tests worth on performing. First trying to give ζ depending from time or distance and then from both. Increasing time and trying more complicated functions are also worth to try. Let's pick f(x) = cos(x) in the range [-π/2 , π/2 ] and start experimenting. Let's fix time for now, say t=15. First is giving time dependent function say ζ = $t^2$. As it appears for t=15 tp2 doesn't give interesting result(*Figure 4.2.4)*, however for t=100 distortion goes intense(*Figure4.2.5*).

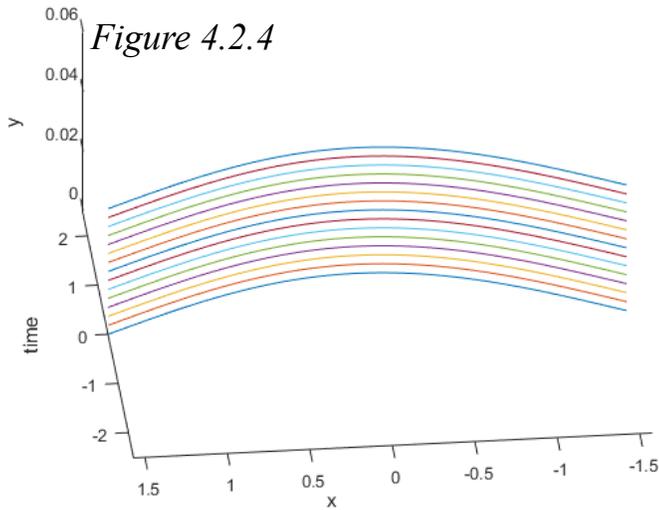

*Figure 4.2.4*

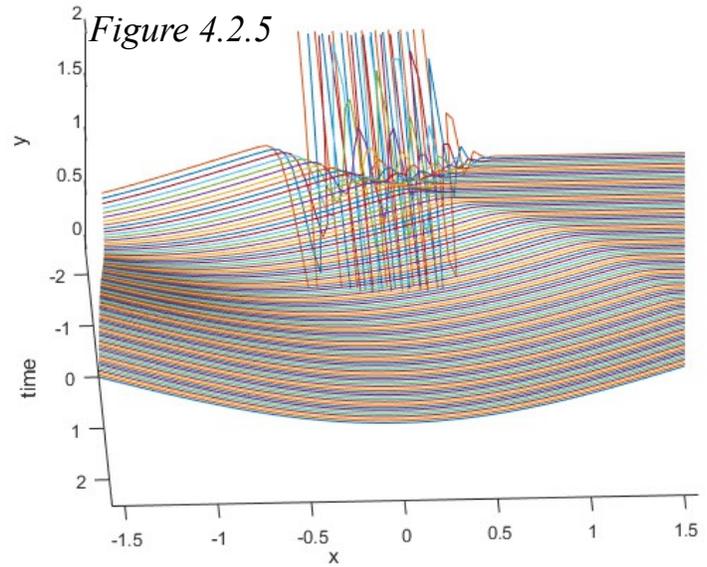

*Figure 4.2.5*

Now let's run do the same using tp3. There is three plots presented, first illustrates the beginning, second is the same wave after some time, and last after long time. AS we see Wave eventually becomes a line which is normal for the wave in nature, however if we wait long enough unnatural behavior is witnessed. (*Figures 4.2.6 a b c*).

*Figures 4.2.6 a b c*

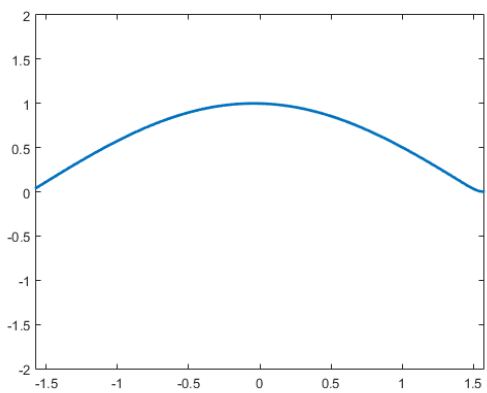 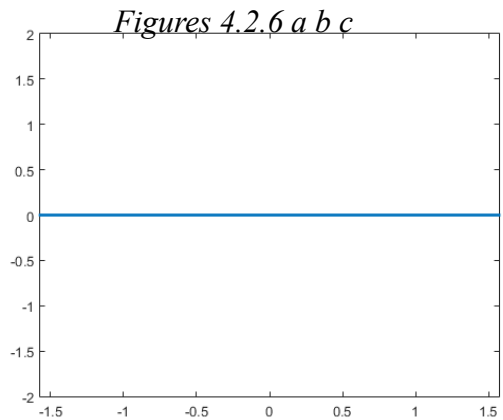 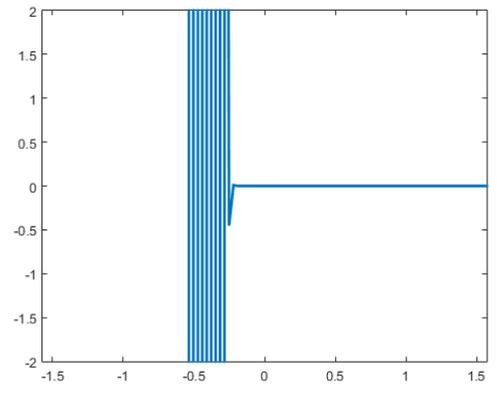

Next Let's take ζ = $x^2$. Also the initial t is 25. In *Figure 4.2.7* we already see the beginning of graininess, However for t = 100, we see that during propagation graininess is gone. And we have a simple crushing wave. (*Figure 4.2.8*). Also constructing by tp3 the beginning is similar with *Figure 4.2.6 a.* The continuation is different (*Figure 4.2.9*) and explains the previous figure.

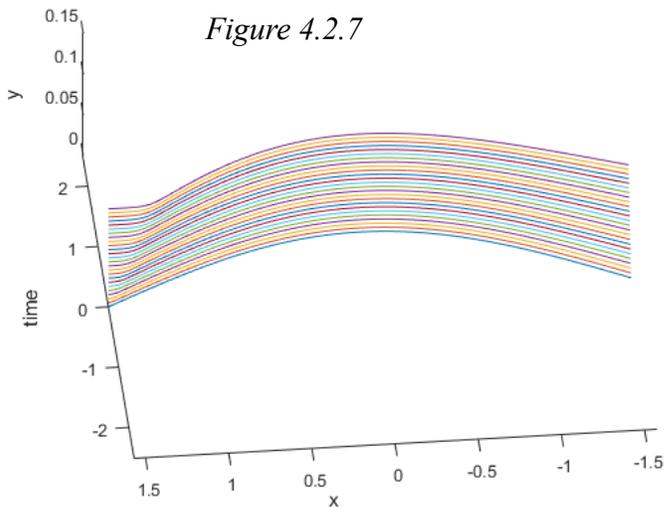

Figure 4.2.7

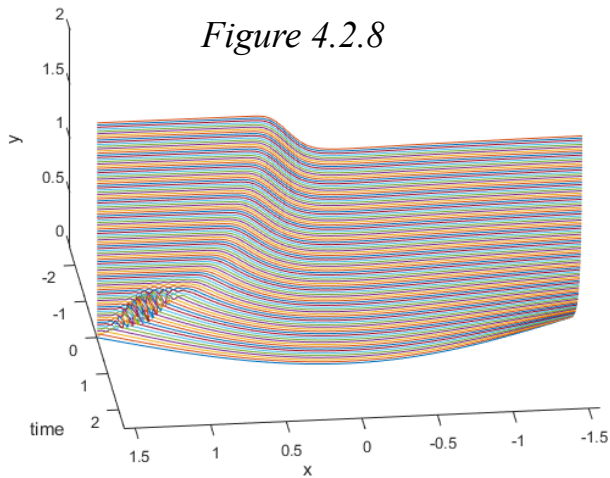

Figure 4.2.8

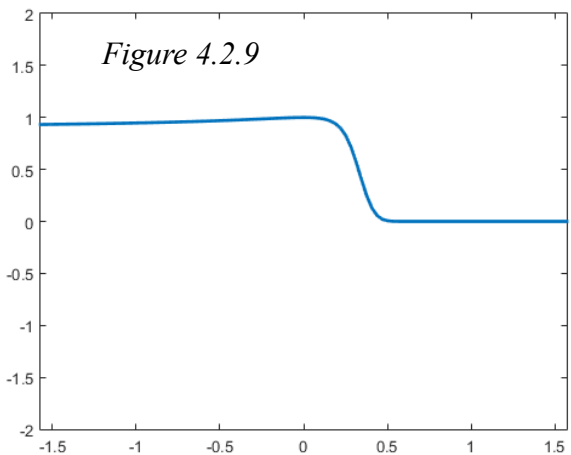

Figure 4.2.9

Next let's pick $\zeta = x^2+t^2$. And perform similar actions.

The examination shows that for tp3 the behavior is similar as in figures *Figures 4.2.6*.

For tp2 the result for t = 25 is similar as *Figure 4.2.7* However for t=100, wave shows hysterical behavior. *Figure 4.2.10*. (Without blue color it was too hard to see the wave)

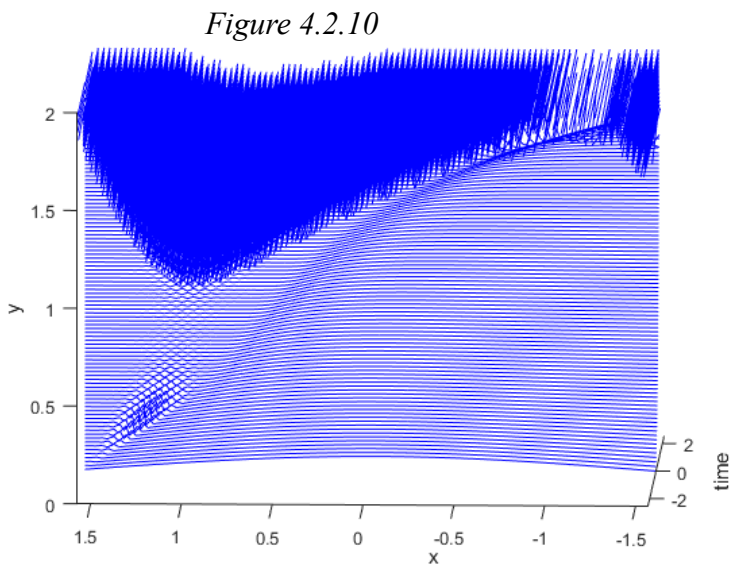

Figure 4.2.10

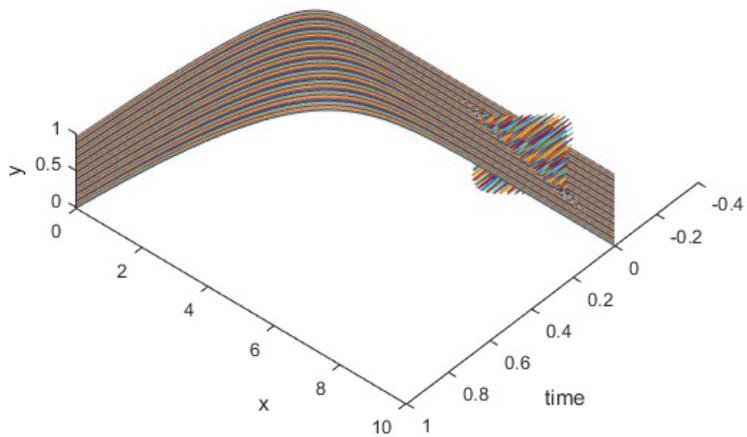

## 4.3 Conclusion

From Figures illustrated in Section 4.2 can be concluded that even though some times the numeric solution for partial derivative equations can give a correct result, it is not always describing a physical phenomena. As Waves simply cannot behave in nature as for example in *figure 4.2.5*.

Few Things can be noticed. First there are cases when error occurs, but later it is self fixed and the continuation is normal.( For instance *Figure 4.2.11* Describing $e^x$ with $\zeta = x$). Another thing that can be noticed is that Waves are inclined to the right, which is mainly caused by our method, double use of the same point).

There are ways which can make this solution more realistic. That is making small $\Delta t$ and/or $\Delta x$. However it will result in slower calculation. There are known ways for stable Calculations. One of such methods is presented by Suren Khachatryan in the paper *"Modification of the Method of Integration along Characteristics as Adaptive Mesh Approach in Solution of Hyperbolic and Parabolic PDEs"*.

# Chapter 5

# Appendices

## 5.1 Glossary

**discretization** – method which involves dividing to small pieces in order to solve.
**shock waves** – abrupt jumps in density, velocity and pressure
**hyperbolic wave** – wave mathematically formulated through hyperbolic partial differential equations is called hyperbolic wave.
**Newton's second law** is F = ma. It says Force(F) on an objects equal to mass(m) multiplied by acceleration (**a**)
**diffusing viscous term** – describes the physical transer quantity through a small area with normal velocity per some time.
**viscous flow** – Viscosity describes fluids resistance to deformation caused by stress (can be shear stress or tensile stress). Equivalent term for metals would be thickness
**shear stress** is given by a formula F/A, where F is force and A is cross sectional area of material with area parallel to the applied force vector. It is applied to fluids which move along solid boundary and is given by formula $\mu \frac{\partial u}{\partial y}$, where $\mu$ is dynamic viscosity, u – velocity along the boundary, y – height above boundary
**tensile stress** – in fluid tensile stress describes the internal force applied on a given particle, caused by neighboring particles.

## 5.2 Our Code

```matlab
function [  ] = OneD(f, c0, t, a, b)
n = 100;
step = (b-a) / n;
k = 0;

while(k<=t)
    x = [a+k zeros(1,n)];
    y=[f(a+k-k*c0) zeros(1,n)];

    for i=1:n
        y(i+1)= f(x(i)-k*c0) ;
        x(i+1)=x(i)+step;
    end
    pause(0.000000002)
    k1=[k k*ones(1,n)];
    plot3(x,k1,y,'b')
    xlabel('x'); ylabel('time'); zlabel('y');
    hold on
    k=k+0.01;
end

end
```

```matlab
function y = tp1(f,g,p1,p2,rt)
nx = 100; %We divide the interval [p1,p2] into nx parts
nt = 5000; %We divide the time interval [0,rt] into nt parts
dx = (p2-p1) / nx; %the length of x-subintervals, x-step
dt = rt / nt; %the length of t-subintervals, t-step
x = p1:dx:p2; %the x-axis division points
t = 0:dt:rt; %the t-axis division points

%Initialization
for m = 1:nx+1
    a(1,m) = f(p1+(m-1)*dx);% filling the first row, the initial data
end

%Calculation
for j = 1:nt
    for i = 2:nx
        a(j+1,i) = dt*g(p1+(i-1)*dx,(j-1)*dt)*(a(j,i+1)-a(j,i-1))/(2*dx) + a(j,i); %filling array by 3points
    end
        a(j+1,1) = a(j+1,2); %not the best way to calculate at the endpoints
        a(j+1,nx) = a(j+1,nx-1);%not the best way to calculate at the endpoints
end
a(i,:);

for i =1:rt
    axis([p1 p2 -2.5 2.5])   %To fix the axis lengths
    t1=[t(i) t(i)*ones(1,nx)];
    pause(0.002);
    plot3(x,a(i,:),t1);
    xlabel('x'); ylabel('time'); zlabel('y');
    hold on
end

end
```

```matlab
function y = tp2(f,g,p1,p2,rt)
nx = 100; %We divide the interval [p1,p2] into nx parts
nt = 5000; %We divide the time interval [0,rt] into nt parts
dx = (p2-p1)/nx; %the length of x-subintervals, x-step
dt = rt/nt; %the length of t-subintervals, t-step
x = p1:dx:p2; %the x-axis division points
t = 0:dt:rt; %the t-axis division points

%Initialization
for m = 1:nx+1
    a(1,m) = f(p1+(m-1)*dx);% filling the first row, the initial data
end

%Calculation
for j = 1:nt
    for i = 1:nx
        a(j+1,i) = dt*g(p1+(i-1)*dx,(j-1)*dt)*(a(j,i+1)-a(j,i))/(dx) + a(j,i); %filling array by 3points
    end
end

for i =1:rt
    t1=[t(i) t(i)*ones(1,nx)];
    pause(0.002);
    plot3(x,a(i,:),t1);
    xlabel('x'); ylabel('time'); zlabel('y');
   hold on
end

end
```

```matlab
function y = tp3(f,g,p1,p2,rt)
nx = 100; %We divide the interval [p1,p2] into nx parts
nt = 5000; %We divide the time interval [0,rt] into nt parts
dx = (p2-p1)/nx; %the length of x-subintervals, x-step
dt = rt/nt; %the length of t-subintervals, t-step
x = p1:dx:p2; %the x-axis division points
t = 0:dt:rt; %the t-axis division points

%Initialization
for m = 1:nx+1
    a(1,m) = f(p1+(m-1)*dx);% filling the first row, the initial data
end

%Calculation
for j = 1:nt
    for i = 1:nx
        a(j+1,i) = dt*g(p1+(i-1)*dx,(j-1)*dt)*(a(j,i+1)-a(j,i))/(dx) + a(j,i); %filling array by 3points
    end
end

%Plotting now
for i =1:nt
  pause(0.0004);
  plot(x,a(i,:),'LineWidth',2);
  axis([p1 p2 -2 2]) %To fix the axis lengths
end

end
```

```matlab
function y = tp4(f,g,p1,p2,rt)
nx = 100; %We divide the interval [p1,p2] into nx parts
nt = 5000; %We divide the time interval [0,rt] into nt parts
dx = (p2-p1)/nx; %the length of x-subintervals, x-step
dt = rt/nt; %the length of t-subintervals, t-step
x = p1:dx:p2; %the x-axis division points
t = 0:dt:rt; %the t-axis division points

%Initialization
for m = 1:nx+1
    a(1,m) =f(p1+(m-1)*dx);% filling the first row, the initial data
end

%Calculation
for j = 1:nt
    for i = 2:nx
        a(j+1,i) = dt*g(p1+(i-1)*dx,(j-1)*dt)*(a(j,i+1)-a(j,i-1))/(2*dx) + a(j,i); %filling array by 3points
    end
    a(j+1,1) = a(j+1,2); %not the best way to calculate at the endpoints
    a(j+1,nx) = a(j+1,nx-1);%not the best way to calculate at the endpoints
end

%Plotting now
for i =1:nt
    pause(0.0004);
    plot(x,a(i,:),'LineWidth',2);
    axis([p1 p2 -2 2]) %To fix the axis lengths
end

end
```

*Published with MATLAB® R2015a*

# Special thanks

Special thanks to our professors:

Special thanks to Michael Poghosyan for his patience and help throughout our project.
Special thanks to Suren Khachatryan for explaining complicated physics behind.